\documentclass[12pt]{article}
\usepackage[margin=2cm]{geometry}
\usepackage{amsfonts}
\usepackage{amssymb}
\usepackage{amsthm}
\usepackage{amsmath}
\usepackage{latexsym}
\usepackage{color}
\usepackage{mathrsfs}

\begin{document}
\allowdisplaybreaks[4]
\newtheorem{lemma}{Lemma}
\newtheorem{pron}{Proposition}
\newtheorem{re}{Remark}
\newtheorem{thm}{Theorem}
\newtheorem{Corol}{Corollary}
\newtheorem{exam}{Example}
\newtheorem{defin}{Definition}
\newtheorem{remark}{Remark}
\newtheorem{property}{Property}
\newcommand{\la}{\frac{1}{\lambda}}
\newcommand{\sectemul}{\arabic{section}}
\renewcommand{\theequation}{\sectemul.\arabic{equation}}
\renewcommand{\thepron}{\sectemul.\arabic{pron}}
\renewcommand{\thelemma}{\sectemul.\arabic{lemma}}
\renewcommand{\there}{\sectemul.\arabic{re}}
\renewcommand{\thethm}{\sectemul.\arabic{thm}}
\renewcommand{\theCorol}{\sectemul.\arabic{Corol}}
\renewcommand{\theexam}{\sectemul.\arabic{exam}}
\renewcommand{\thedefin}{\sectemul.\arabic{defin}}
\renewcommand{\theremark}{\sectemul.\arabic{remark}}
\def\REF#1{\par\hangindent\parindent\indent\llap{#1\enspace}\ignorespaces}
\def\lo{\left}
\def\ro{\right}
\def\be{\begin{equation}}
\def\ee{\end{equation}}
\def\beq{\begin{eqnarray*}}
\def\eeq{\end{eqnarray*}}
\def\bea{\begin{eqnarray}}
\def\eea{\end{eqnarray}}
\def\d{\Delta_T}
\def\r{random walk}
\def\o{\overline}

\title{\large\bf Some positive conclusions related to the Embrechts-Goldie conjecture}
\author{\small Yuebao Wang$^1$\thanks{Research supported by National Natural Science Foundation of China
(No.s 11071182).}
\thanks{Corresponding author.
Telephone: +86 512 67422726. Fax: +86 512 65112637. E-mail:
ybwang@suda.edu.cn}~~Zhaolei Cui$^2$~~Hui Xu$^1$
\\
{\footnotesize\it 1 School of Mathematical Sciences, Soochow
University, Suzhou 215006, China}\\
{\footnotesize\it 2 School of mathematics and statistics, Changshu Institute of Technology, Suzhou 215000, China}}
\date{}
\maketitle

\begin{center}
{\noindent\small {\bf Abstract }}
\end{center}

{\small
In this paper, we give some conditions, under which, if an infinitely divisible distribution supported on $[0,\infty)$ belongs to the intersection of exponential distribution class $\mathcal{L}(\gamma)$ for some $\gamma\ge0$ and generalised subexponential distribution class $\mathcal{OS}$, then its L$\rm\acute{e}$vy spectral distribution or convolution of the distribution with itself also belongs to the same one. To this end, we discuss the closure under the compound convolution roots for the class. In addition, we do some in-depth discussion about the above-mentioned conditions, and provide some types of distributions satisfying them. Further, we obtain some local versions of the above-mentioned results by the Esscher transform of distributions. Therefore, some positive conclusions related to the Embrechts-Goldie conjecture are obtained. Prior to this, all corresponding results are negative.

\medskip

\medskip

}

\section{Introduction}
\setcounter{thm}{0}\setcounter{Corol}{0}\setcounter{lemma}{0}\setcounter{pron}{0}\setcounter{equation}{0}
\setcounter{remark}{0}\setcounter{exam}{0}\setcounter{property}{0}\setcounter{defin}{0}

Let $H$ be an infinitely divisible distribution supported on $[0,\infty)$
with the Laplace transform
\begin{eqnarray}\label{500}
\int_0^\infty \exp\{-\lambda y\}H(dy)=\exp\Big\{-a\lambda-\int_0^\infty(1-e^{\lambda y})\upsilon(dy)\Big\},
\end{eqnarray}
where $a\ge0$ is a constant, and $\upsilon$ is a Borel measurable function supported on $(0,\infty)$ with the properties $\upsilon((1,\infty))<\infty$ and $\int_0^1y\upsilon(dy)<\infty$, which is called the L$\rm\acute{e}$vy spectral measure. Let
$$F(x)=\upsilon(x)\textbf{1}(x>1)/\mu=\upsilon(0,x]\textbf{1}(x>1)/\mu$$
for all $x\in(-\infty,\infty)$ be the L$\rm\acute{e}$vy spectral distribution
generated by the measure $\upsilon$.
The distribution $H$ admits the representation $H=H_1*H_2$, which is reserved
for convolution of two distributions $H_1$ and $H_2$ 
satisfying $\overline H_1(x)=O(e^{-\beta x})$ for some constant $\beta>0$ and
\begin{eqnarray}\label{501}
H_2(x)=e^{-\mu}\sum_{k=0}^\infty F^{*k}(x)\mu^k/k!
\end{eqnarray}
for all $x\in(-\infty,\infty)$, where $\mu=\upsilon\big((1,\infty)\big)$, $F^{*k}$
is the $k$-fold convolution of $F$ with itself for all integers $k\ge1$ and $F^{*0}$
is the distribution degenerate at zero. See, for example, Feller \cite{F1971}.

We might also say that $F$ is an ``input" and that $H$ is an ``output" in a system. Usually, we use ``input" $F$ to infer the ``output" $H$. However, when the $F$ is in a ``black box", then we need to use $H$ to infer $F$. In this paper, our main research topic is that, under what conditions, a L$\rm\acute{e}$vy spectral distribution or convolution of the distribution with self belongs to certain distribution class, if the corresponding infinite divisible distribution belongs to the same one? In this way, we first need to recall the concepts and notations of some distribution classes.

Here and later, without special statement, all limits are taken as $x$ tends to infinity. 
And for a distribution $V$, let $\overline{V}=1-V$ be tailed distribution of $V$.

For some constant $\gamma\ge0$, a distribution $V$ supported on $[0,\infty)$ or $(-\infty,\infty)$ belongs to the distribution class $\mathcal{L}(\gamma)$, if $\overline V(x)>0$ for all $x$ and $\lim\overline V(x-t)/\overline V(x)=e^{\gamma t}$ denoted by
$$\overline V(x-t)\sim \overline V(x)e^{\gamma t}$$
for any fixed $t>0$.

In the above definition, if $\gamma>0$ and the distribution $V$ is lattice, then $x$ and $t$ should be restricted to values of the lattice span, see Bertoin and Doney \cite{BD1996}. As everyone knows that, if $V\in\mathcal{L}(\gamma)$, then
$$\mathcal{H}(V,\gamma)=\{h(\cdot):h(x)\uparrow\infty,h(x)/x\to0,\overline V(x-t)\sim \overline V(x)e^{\gamma t}\ \text{uniformly for all}\ |t|\le h(x)\}\neq\phi.$$

A distribution $V$ supported on $[0,\infty)$ or $(-\infty,\infty)$ belongs to the distribution class $\mathcal{S}(\gamma)$ for some $\gamma\ge0$,
if $M(V,\gamma)=\int_0^{\infty}e^{\gamma y}V(dy)\ or\ \int_{-\infty}^{\infty}e^{\gamma y}V(dy)$ is finite, $V\in\mathcal{L}(\gamma)$ and
$$\overline {V^{*2}}(x)\sim 2M(V,\gamma)\overline V(x).$$

In some literatures, these two classes are called exponential distribution class and convolution equivalent distribution class, respectively. In particular, the classes $\mathcal{L}(0)$ and $\mathcal{S}(0)$ are called long-tailed distribution class and subexponeantial distribution class, denoted by $\mathcal{L}$ and $\mathcal{S}$, respectively. It should be noted that the requirement $V\in\mathcal{L}$ is not needed in the definition of the class $\mathcal{S}$ when $F$ is supported on $[0,\infty)$.

The class $\mathcal{S}$ was introduced by Chistyakov \cite{C1964} and the class $\mathcal{S}(\gamma)$ for some $\gamma>0$ by Chover et al. \cite{CNW1973a,CNW1973b} for the support $[0,\infty)$, and Tang and Tsitsiashvili \cite{TT2003} or Pakes \cite{P2004} for the support $(-\infty,\infty)$. The classes $\cup_{\gamma\ge0}\mathcal{S}(\gamma)$ and $\cup_{\gamma\ge0}\mathcal{L}(\gamma)$ are properly included in the following two distribution classes, respectively.

A distribution $V$ supported on $[0,\infty)$ or $(-\infty,\infty)$ belongs to the generalised subexponential distribution class
$\mathcal{OS}$ introduced by Kl\"{u}ppelberg \cite{K1990} or Shimura and Watanabe \cite{SW2005}, if
$$C^*(V)=\limsup\overline{V^{*2}}(x)/\overline{V}(x)<\infty.$$

A distribution $V$ supported on $[0,\infty)$ or $(-\infty,\infty)$ belongs to the generalised long-tailed distribution class $\mathcal{OL}$ introduced by Shimura and Watanabe \cite{SW2005}, if for any positive constant $t$,
$$C^*(V,t)=\limsup\overline{V}(x-t)/\overline{V}(x)<\infty.$$

Further, Shimura and Watanabe \cite{SW2005} show that the inclusion relation $\mathcal{OS}\subset\mathcal{OL}$ is proper.

Recall the research of the above-mentioned our main research topic, for the class $\mathcal{S}(\gamma)$, Embrechts et al. \cite{EGV1979} for $\gamma=0$, Sgibnev \cite{S1990} and Pakes \cite{P2004} for $\gamma>0$ have already had positive result which show that, the L$\rm\acute{e}$vy spectral distribution $F$ of an infinitely divisible distribution $H$ belongs to the class $\mathcal{S}(\gamma)$ when $H\in\mathcal{S}(\gamma)$ combined with some conditions. For other distribution classes, however, there are only some negative results that, here exists an infinitely divisible distribution $H$ which belongs to the class, while its L$\acute{e}$vy spectral distribution $F$ is not, see Theorem 1.1 (iii) of Shimura and Watanabe \cite{SW2005} for class $\mathcal{OS}$, Theorem 1.2 (3) of Xu et al. \cite{XFW2015} for class $\mathcal{L}\cap\mathcal{OS}$ and Theorem 1.1 of Xu et al. \cite{XWCY2016} for class $\mathcal{L}(\gamma)\cap\mathcal{OS}$ with some $\gamma>0$. Therefore, for the class $\mathcal{L}(\gamma)\cap\mathcal{OS}$, more precisely, for the class $\mathcal{L}(\gamma)\cap\mathcal{OS}\setminus\mathcal{S}(\gamma)$, the following question is raised naturally:\\

\textbf{Problem 1.1.} Under what conditions, the L$\rm\acute{e}$vy spectral distribution or convolution of the distribution with self belongs to the class $\mathcal{L}(\gamma)\cap\mathcal{OS}$ for some $\gamma\ge0$, if corresponding infinitely divisible distribution $H$ belongs to the same one?\\

For this problem, we give a positive answer as follows. 

\begin{thm}\label{thm40}
Let $H$ be an infinitely divisible distribution supported on $[0,\infty)$
with the Laplace transform (\ref{500}) and L$\acute{e}$vy spectral distribution $F$. Assume that $H\in{\mathcal{L}}(\gamma)\cap\mathcal{OS}$ for some $\gamma\ge0$ and $\overline{H_1}(x)=o\big(\overline{H_2}(x)\big)$. And then suppose that, for all $k\ge1$,
\begin{eqnarray}\label{thm102}
\liminf\overline{F^{*k}}(x-t)/\overline{F^{*k}}(x)\ge e^{\gamma t}\ for\ all\ t>0.
\end{eqnarray}
Then $H_2\in{\mathcal{L}}(\gamma)\cap\mathcal{OS}$ and $\overline{H_2}(x)\thickapprox\overline{H}(x)$. And there exists a integer $n_0\ge1$ such that $F^{*n}\in{\mathcal{L}}(\gamma)\cap\mathcal{OS}$ for all $n\ge n_0$.

Further, if $n_0\ge 2$ and there exits an integer $1\le l_0\le n_0-1$ such that $F^{*l_0}\in\mathcal{OS}$, then $F^{*n}\in{\mathcal{L}}(\gamma)\cap\mathcal{OS}$ for all $n\ge l_0$.

In particular, if $l_0=1$, that is $F\in\mathcal{OS}$, then $F\in\cal{L}(\gamma)\big(\cap\mathcal{OS}\big)$.
\end{thm}

\begin{remark}\label{remark1010}
i) Clearly, when $\gamma=0$, (\ref{thm102}) holds automatically for all $k\ge1$. In Theorem 2.2 of Xu et al. \cite{XFW2015}, there is a distribution $F$ satisfying (\ref{thm102}) such that,
$H_2$ and $H$ belong to the class $({\mathcal{L}}\cap\mathcal{OS})\setminus\mathcal{S}$,
while $F\in\mathcal{OL}\setminus({\mathcal{L}}\cap\mathcal{OS})$. Because $F^{*2}\in\mathcal{OS}$, that is $l_0=2$, Theorem \ref{thm40} implies that $F^{*n}\in\mathcal{L}\cap\mathcal{OS}$ for all $n\ge l_0=2$. A similar example can be found in Proposition 2.1 of Xu et al. \cite{XFW2015}, where $F\notin\mathcal{OL}$ with $l_0=2$.

When $\gamma>0$, the condition (\ref{thm102}) for $k=1$ has been used in Lemma 7 and Theorem 7 of Foss and Korshunov \cite{FK2007}. And Proposition 4.3 of Xu et al. \cite{XWCY2016} provides a type of distribution $F$ satisfying (\ref{thm102}) such that, $H_2$ and $H$ belong to the class $\big(\mathcal{L}(\gamma )\cap\mathcal{OS}\big)\setminus\mathcal{S}(\gamma)$, and $F\in\cal{OL}\setminus\big(\mathcal{L}(\gamma)\cup\mathcal{OS}\big)$, while $F^{*n}\in\mathcal{L}(\gamma )\cap\mathcal{OS}$ for all $n\ge l_0=2$ by Theorem \ref{thm40}.


Further discussion on the condition (\ref{thm102}) is put in Subsection 6.1 of present paper. There, we try to replace the condition (\ref{thm102}) for all integers $k\ge1$ with more simple and concrete condition.

ii) In the above-mentioned examples, $F\notin\mathcal{L}(\gamma)\cup\mathcal{OS}$ for some $\gamma\ge0$. Therefore, under the condition (\ref{thm102}), the condition that $F\in\mathcal{OS}$ is necessary in certain sense for $F\in\cal{L}(\gamma)$.

iii) In this theorem, the main object of research is the class $\mathcal{L}(\gamma)\cap\mathcal{OS}$ for some $\gamma\ge0$. We note that, here many distributions in the class $(\mathcal{L}(\gamma)\cap\mathcal{OS})\setminus\mathcal{S}(\gamma)$ have been found, see, for example, Leslie \cite{L1989}, Kl\"{u}ppelberg and Villasenor \cite{KV1991}, Shimura and Watanabe \cite{SW2005}, Lin and Wang \cite{LW2012}, Wang et al. \cite{WXCY2016}, Xu et al. \cite{XFW2015} and Xu et al. \cite{XWCY2016}. For research on $\mathcal{OS}$, besides the above-mentioned literatures, the reader can refer to Watanabe and Yamamura \cite{WY2010}, Yu and Wang \cite{YW2014}, Beck et al. \cite{BBS2015}, Xu et al. \cite{XSW2015, XSWC2016}, etc.
\end{remark}

From (\ref{501}) we can find that, in order to prove Theorem \ref{thm40}, we first need to solve the following Problem 1.2 involving compound distribution or compound convolution. Let $\tau$ be a nonnegative integer-valued random variable with masses $p_k=\textbf{P}(\tau=k)$ for all integers $k\ge0$ satisfying $\sum_{k=0}^\infty p_k=1$. And let $V$ be a distribution. Write a compound convolution generated by the random variable and distribution
\begin{eqnarray}\label{101}
V^{*\tau}=\sum_{k=0}^\infty p_kV^{*k}.
\end{eqnarray}
For convenience, in this paper, we set up $p_k>0$ for all integers $k\ge0$. In fact, if $\tau$ is a nonnegative integer-valued random variable with masses $p_{k_m}>0$ for all integers $m\ge1$ satisfying $\sum_{m=0}^\infty p_{k_m}=1$, where $k_1=1$, then all conclusions of the paper still hold.


\textbf{Problem 1.2.} Under what conditions, 
the 
distribution $V$ or its convolution with self belongs to the class $\mathcal{L}(\gamma)\cap\mathcal{OS}$, if $V^{*\tau}\in\mathcal{L}(\gamma)\cap\mathcal{OS}$?

Generally, it is called a topic on closure under compound convolution roots for some distribution class.

It is well known that, compound convolution including its convolution with other distribution has extensive and important applications in various fields, such as risk model, queuing system, branching process, infinitely divisible distribution, and so on. See, for example, Embrechts et al. \cite{EKM1997} and Foss et al. \cite{FKZ2013}. 

The topic in Problem 1.2 is a natural extension on the famous conjecture in Embrechts and Goldie \cite{EG1980,EG1982} for the distribution class $\mathcal{L}(\gamma)$ for some $\gamma\ge0$. Some of the latest results about the conjecture and the related problems can be found in Xu et al. \cite{XFW2015}, Watanabe  \cite{W2015}, Watanabe and Yamamuro \cite{WY2016}, Xu et al. \cite{XWCY2016}.

In the references mentioned above, Theorem 2.2 of Xu et al. \cite{XFW2015} for $\gamma=0$ and Proposition 4.3 of Xu et al. \cite{XWCY2016} for $\gamma>0$ show that,
the distribution class $\mathcal{L}(\gamma)\cap\mathcal{OS}$ is not closed under compound convolution roots. These conclusions give some negative answers to the Embrechts-Goldie conjecture and Problem 1.2.

Finally, we naturally hope to answer the following question.

\textbf{Problem 1.3.} Will there be some results similarly to Theorem \ref{thm40} for some local distribution classes?

In Section 3, we prove Theorem \ref{thm40}. To this end, we give a positive answer 
for Problem 1.2 in Section 2. And in Section 4, by the Esscher transform between distributions and the related results of Wang and Wang \cite{WW2011}, we get three positive results for Problem 1.3, 
their proofs are in Section 5. Finally, in Section 6, we respectively provide some more specific and convenient conditions which can replace the conditions (\ref{thm102}) in Theorem \ref{thm40} and (\ref{thm103}) in Theorem \ref{thm1} below.

\section{On the compound convolution} 
\setcounter{thm}{0}\setcounter{Corol}{0}\setcounter{lemma}{0}\setcounter{pron}{0}\setcounter{equation}{0}
\setcounter{remark}{0}\setcounter{exam}{0}\setcounter{property}{0}\setcounter{defin}{0}

In the following, it is agreed that all distributions are supported on $[0,\infty)$. Let $V$ be a distribution. And random variable $\tau$ and compound convolution $V^{*\tau}$ are as defined in Section 1.

Now, we give a positive answer for Problem 1.2.

\begin{thm}\label{thm1}
For each constant $0<\varepsilon<1$, assume that there is an integer $n_0=n_0(V,\tau,\varepsilon)\ge 1$ such that,
\begin{eqnarray}\label{thm103}
\sum\limits_{k=n_0+1}^{\infty}p_{k}\overline{V^{*(k-1)}}(x)\leq\varepsilon\overline{V^{*\tau}}(x)\ for\ all\ x\ge0.
\end{eqnarray}
If $V^{*\tau}\in\mathcal{OS}$, then $V^{*n}\in\mathcal{OS}$ for all $n\ge n_0$.

Moreover, suppose that (\ref{thm102}) holds for all $k\ge1$. If $V^{*\tau}\in{\mathcal{L}}(\gamma)$, then for all $n\ge n_0$, $V^{*n}\in{\mathcal{L}}(\gamma)$, thus $V^{*n}\in{\mathcal{L}}(\gamma)\cap\mathcal{OS}$. 

Here, if $n_0\ge2$ and if there exists an integer $1\le l_0\le n_0-1$ such that $V^{*l_0}\in\cal{OS}$, then $V^{*n}\in\cal{L}(\gamma)\cap\mathcal{OS}$ for all $n\ge l_0$.

In particular, if $V\in\mathcal{OS}$, that is $l_0=1$, then $V\in\cal{L}(\gamma)(\cap\mathcal{OS})$.
\end{thm}
\begin{remark}\label{remark101}
i) The condition (\ref{thm103}) was used in Watanabe and Yamamuro \cite{WY2010}, Yu and Wang \cite{YW2014}, Xu et al. \cite{XSW2015} and Xu et al. \cite{XWCY2016}. As Watanabe and Yamamuro \cite{WY2010} points out that, if $p_{k+1}/p_k\to0$ as $k\to\infty$, for example, $p_k=e^{-\lambda}\lambda^k/k!$ for $k\ge0$, then condition (\ref{thm103}) is satisfied. More distributions satisfying the condition (\ref{thm103}), for which the condition that $p_{k+1}/p_k\to0$ as $k\to\infty$ is fail, can be found in Subsection 6.2 below.

ii) Of course, we want to find a minimal $n_0$. Clearly, if $\varepsilon$ becomes larger, then $n_0$ becomes smaller. In Remark \ref{remark1010}, however, we find $n_0=l_0=2$ for some compound convolution. In other words, we cannot guarantee $n_0=1$ for each compound convolution.
\end{remark}

\proof In order to prove Theorem \ref{thm1}, we first give an equivalent form of the condition (\ref{thm103}) in the case that $F^{*\tau}\in\mathcal{OS}$.

\begin{lemma}\label{lemma1}
If $V^{*\tau}\in\mathcal{OS}$ with $p_1=P(\tau=1)>0$, then the following two propositions are equivalent to each other.

i) For any $0<\varepsilon<1$, there is an integer $n_0=n_0(V,\tau,\varepsilon)\ge1$ such that
\begin{eqnarray}\label{l101}
\sum\limits_{k=n_0+1}^{\infty}p_k\overline{V^{*k}}(x)\leq\varepsilon\overline{V^{*\tau}}(x)\ for\ all\ x\ge0.
\end{eqnarray}

ii) For any $0<\varepsilon<1$, there is an integer $n_0=n_0(F,\tau,\varepsilon)\ge1$ such that (\ref{thm103}) holds.
\end{lemma}
\proof  We only need to prove ii) $\Longrightarrow$ i). To this end, we denote  $$D^*(V^{*\tau})=\sup_{x\ge0}\overline{(V^{*\tau})^{*2}}(x)/\overline{V^{*\tau}}(x).$$
Clearly, $1\le D^*(V^{*\tau})<\infty$ by $V^{*\tau}\in\mathcal{OS}$. For any $0<\varepsilon<1$, there is a $0<\varepsilon_1<p_1/D^*(V^{*\tau})<1$ such that $\varepsilon=\varepsilon_1D^*(V^{*\tau})/p_1$. For above $\varepsilon_1$, by proposition ii), there is an integer $n_0=n_0(V,\tau,\varepsilon_1)\ge1$ such that (\ref{thm103}) holds. Further, we assume that
\begin{eqnarray}\label{l1031}
a_{n_0}=\sum\limits_{k=n_0+1}^{\infty}p_k<\varepsilon_1.
\end{eqnarray}
Define a distribution $G_{n_0}$ by $G_{n_0}(x)=\sum\limits_{k=n_0+1}^\infty p_{k}V^{*(k-1)}(x)/a_{n_0}$ for all $x$. Then by (\ref{thm103}) and (\ref{l1031}), for all $x\ge0$, we have
\begin{eqnarray*}
\sum\limits_{k=n_0+1}^{\infty}p_{k}\overline{V^{*k}}(x)&=&\overline{V*\sum\limits_{n_0+1\le k<\infty}p_{k}V^{*(k-1)}}(x)\nonumber\\
&=&a_{n_0}\int_{0-}^x\overline{G_{n_0}}(x-y)V(dy)+a_{n_0}\overline{V}(x)\nonumber\\
&\le&\varepsilon_1\Big(\int_{0-}^x\overline{V^{*\tau}}(x-y)V(dy)+\overline{V}(x)\Big)\nonumber\\
&=&\varepsilon_1\overline{V*V^{*\tau}}(x)\nonumber\\
&\le&\varepsilon_1 D^*(V^{*\tau})\overline{V^{*\tau}}(x)/p_1
=\varepsilon\overline{V^{*\tau}}(x),
\end{eqnarray*}
that is we get proposition i).\hfill$\Box$\\

Now, we prove the first conclusion of Theorem \ref{thm1}. For some $0<\varepsilon_0<1$, by (\ref{thm103}) with $\varepsilon=\varepsilon_0$, $V^{*\tau}\in\mathcal{OS}$ and Lemma \ref{lemma1}, there is an integer $n_0=n_0(V,\tau,\varepsilon_0)\ge1$ such that
$$\sum\limits_{k=1}^{n_0}p_k\overline{V^{*k}}(x)\geq(1-\varepsilon_0)\overline{V^{*\tau}}(x)\ \text{for all}\ x\ge0,$$
thus $\overline{V^{*\tau}}(x)\approx\overline{V^{*n_0}}(x).$
Then by $V^{*\tau}\in\cal{OS}$, we immediately get $V^{*n_0}\in\cal{OS}$. According to Proposition 2.6 of Shimura and Watanabe (2005) and $V^{*n_0}\in\mathcal{OS}$, $V^{*n}\in\mathcal{OS}$ and $\overline{V^{*\tau}}(x)\approx\overline{V^{*n}}(x)$ for all $n\ge n_0$.

Next, by Lemma \ref{lemma1} and (\ref{thm103}), for any $0<\varepsilon<\varepsilon_0$ and any fixed $n\ge n_0$, there exists an integer $m_0=m_0(V,\tau,\varepsilon)\ge n$ such that $$\sum\limits_{k=m_0+1}^{\infty}p_k\overline{V^{*k}}(x)\leq\varepsilon\overline{V^{*\tau}}(x)\ \text{for all}\ x\ge0.$$
Further, by $V^{*\tau}\in\cal{L}(\gamma)$ and (\ref{thm102}), for the above $0<\varepsilon<\varepsilon_0$ and any fixed $t>0$, there is a constant $x_0=x_0(F,\tau,\varepsilon,t)$ such that, for all $x>x_0$,
\begin{eqnarray*}
&&\varepsilon\overline {V^{*\tau}}(x)\ge\overline {V^{*\tau}}(x-t)-e^{\gamma t}\overline {V^{*\tau}}(x)\nonumber\\
&=&\Big(\sum\limits_{1\le k\neq n\le m_0}+\sum\limits_{k=n}+\sum\limits_{k\ge m_0+1}\Big)p_k\big(\overline{V^{*k}}(x-t)-e^{\gamma t}\overline{F^{*k}}(x)\big)\nonumber\\
&\geq&-\varepsilon e^{\gamma t}\sum\limits_{1\le k\le m_0}p_k\overline{V^{*k}}(x)
+p_{n}\big(\overline{V^{*n}}(x-t)-e^{\gamma t}\overline{V^{*n}}(x)\big)-\varepsilon e^{\gamma t}\overline{V^{*\tau}}(x),
\end{eqnarray*}
which implies that, for all $x>x_{0}$,
\begin{eqnarray*}
\overline{V^{*n}}(x-t)\le e^{\gamma t}\overline{V^{*n}}(x)+(1+2e^{\gamma t})\varepsilon\overline {V^{*\tau}}(x)/p_{n}.
\end{eqnarray*}
Hence by $\overline{F^{*\tau}}(x)\approx\overline{V^{*n}}(x)$, $F^{*n}\in\mathcal{OS}$, $V^{*\tau}\in\mathcal{L}(\gamma)$ and arbitrariness of $\varepsilon$, we can get
\begin{eqnarray}\label{t1041}
\limsup\overline{V^{*{n}}}(x-t)/\overline{V^{*{n}}}(x)\le e^{\gamma t}.
\end{eqnarray}
Thus, combining (\ref{t1041}) and (\ref{thm102}), $V^{*n}\in{\cal{L}}(\gamma)$ for all $n\ge n_0$.

Furthermore, if $V^{*l_0}\in\cal{OS}$ for some $1\le l_0\le n_0-1$, then for any $n\ge l_0$, there are a integer $l_1$ and positive constant $M=M(V,l_0,l_1)$ such that
$l_0\le 2\max\{n_0,n\}\le l_1l_0$ and
$$\overline{V^{*l_0}}(x)\le\overline{V^{*n}}(x)\le\overline{V^{*2n}}(x)\le\overline{V^{*l_0l_1}}(x)\le M\overline{V^{*l_0}}(x).$$
Thus, $V^{*n}\in\cal{OS}$ and $\overline{V^{*n}}(x)\approx\overline{V^{*l_0l_1}}(x)\approx\overline{V^{*\tau}}(x)$. And then use the methods mentioned above, we have $V^{*n}\in{\cal{L}}(\gamma)$ for all $n\ge l_0$.

Particularly, if $l_0=k_0=1$, that is $V\in\cal{OS}$, then $V\in{\cal{L}}(\gamma)$. \hfill$\Box$\\

\section{The proof of Theorem \ref{thm40}}
\setcounter{thm}{0}\setcounter{Corol}{0}\setcounter{lemma}{0}\setcounter{pron}{0}\setcounter{equation}{0}
\setcounter{remark}{0}\setcounter{exam}{0}\setcounter{property}{0}\setcounter{defin}{0}

We first prove the following lemma with a more general form than the first part of Theorem \ref{thm40}, which is a key to the proof of the theorem and has its own independent value.

Let $G_1$ be a distribution. Write $G=G_1*G_2$, where $G_2=V^{*\tau}$ is a compound convolution generated by some distribution $V$ and nonnegative integer-valued random variables $\tau$. Recall that all distributions are supported on $[0,\infty)$.
\begin{lemma}\label{lemma501}
Assume that $G\in{\mathcal{L}}(\gamma)\cap\mathcal{OS}$ for some $\gamma\ge0$ and $\overline{G_1}(x)=o\big(\overline{G_2}(x)\big)$. Further, suppose that conditions (\ref{thm103}) for any $0<\varepsilon<1$ and (\ref{thm102}) for all $k\ge1$ are satisfied. Then $\overline{G_2}(x)\thickapprox\overline{G}(x)$ and $G_2\in{\mathcal{L}}(\gamma)\cap\mathcal{OS}$.
\end{lemma}
\proof First, we prove $\overline{G}(x)\approx\overline{G_2}(x)$. For any $0<\varepsilon<1/\big(2C^*(G)\big)$, by $G\in{\mathcal{L}}(\gamma)$ and $\overline{G_1}(x)=o\big(\overline{G_2}(x)\big)$, there is a constant $A>0$ large enough such that, when $x\ge A$, we have
\begin{eqnarray*}\label{5002}
\overline{G}(x)&=&\int_{0-}^{x-A}\overline{G_1}(x-y)G_2(dy)+\int_{0-}^{A}\overline{G_2}(x-y)G_1(dy)
+\overline{G_1}(A)\overline{G_2}(x-A)\nonumber\\
&\le&\varepsilon\int_{0-}^{x-A}\overline{G}(x-y)G_2(dy)+G_1(A)\overline{G_2}(x-A)
+\overline{G_1}(A)\overline{G_2}(x-A)\nonumber\\
&\le&2\varepsilon C^*(G)\overline{G}(x)+\overline{G_2}(x-A).
\end{eqnarray*}
Therefore, since $G\in\mathcal{OS}\subset\mathcal{OL}$, we have
\begin{eqnarray*}
\big(1-2\varepsilon C^*(G)\big)\overline{G}(x)/C^*(V,A)&\lesssim&\big(1-2\varepsilon C^*(G)\big)\overline{G}(x+A)\le\overline{G_2}(x),
\end{eqnarray*}
that is
$\overline{G}(x)\approx\overline{G_2}(x)$, thus $G_2\in\mathcal{OS}$.

Next, we prove $G_2\in\mathcal{L}(\gamma)$. By (\ref{thm103}), $G_2\in\mathcal{OS}$ and Lemma \ref{lemma1}, there is an integer $n_0\ge k_0$, such that (\ref{l101}) holds. From
(\ref{thm102}), for any $0<\varepsilon<1$ and any $t>0$, there is a constant $x_0=x_0(F,\varepsilon,t)$ such that for all $x\ge x_0$,
\begin{eqnarray}\label{504}
e^{\gamma t}\overline{V^{*k}}(x)/\overline{V^{*k}}(x-t)\le 1+\varepsilon,\ \text{for all}\ k_0\le k\le n_0.
\end{eqnarray}
By (\ref{504}) and (\ref{l101}), we have
\begin{eqnarray*}
\big(e^{\gamma t}\overline{G_2}(x)-\overline{G_2}(x-t)\big)/\overline{G_2}(x-t)
&\le&\sum_{k=k_0}^{n_0} \Big(\frac{e^{\gamma t}\overline{V^{*k}}(x)}{\overline{V^{*k}}(x-t)}-1\Big)+(e^{\gamma t}+1/p_{k_0})\varepsilon\nonumber\\
&\leq&(n_0+e^{\gamma t}+1/p_{k_0})\varepsilon
\end{eqnarray*}
which implies
\begin{eqnarray}\label{506}
\limsup(e^{\gamma t}\overline{G_2}(x)-\overline{G_2}(x-t))/\overline{G_2}(x-t)\le0.
\end{eqnarray}

On the other hand, for any $t>0$, from (\ref{506}), $\overline{G_1}(x)=o\big(\overline{G_2}(x)\big)$ and $G\in\mathcal{L}(\gamma)\cap\mathcal{OS}$, there is a enough large constant $B>2t$   such that, when $x\ge 3B$, we have
\begin{eqnarray}\label{502}
&&e^{2\gamma t}\overline{G}(x)-\overline{G}(x-2t)\leq\Big(\int_{0}^{x-B}+\int_{x-B}^{x}\Big)e^{2\gamma t}\overline{G_1}(x-y)G_2(dy)+e^{2\gamma t}\overline{G_2}(x)-\overline{G_2}(x-2t)\nonumber\\
& &\ \ \ \ \ \ \ \ \ \ \ -\int_{x-2t-B}^{x-2t}\overline{G_1}(x-2t-y)G_2(dy)\nonumber\\
&\leq&\varepsilon e^{2\gamma t}\overline{G^{*2}}(x)+\int_{0}^{B}\Big(e^{2\gamma t}\overline{G_2}(x-y)-\overline{G_2}(x-2t-y)\Big)G_1(dy)\nonumber\\
& &\ \ \ \ \ \ \ \ \ \ \ +\overline{G_1}(B)\big(e^{2\gamma t}\overline{G_2}(x-B)-\overline{G_2}(x-2t-B)\big)\nonumber\\
&\leq&2\varepsilon(1+\varepsilon)e^{2\gamma t}C^*(G)\overline{G}(x)+\Big(\int_{0}^{t}+\int_{t}^{2t}+\int_{2t}^{B}\Big)\Big(e^{2\gamma t}\overline{G_2}(x-y)-\overline{G_2}(x-2t-y)\Big)G_1(dy)\nonumber\\
&\leq&3\varepsilon(1+\varepsilon)e^{2\gamma t}C^*(G)\overline{G}(x)+\int_{t}^{2t}\Big(e^{2\gamma t}\overline{G_2}(x-y)-e^{\gamma y}\overline{G_2}(x-2t)\Big)G_1(dy)\nonumber\\
& &\ \ \ \ \ \ \ \ \ \ \ +\int_{t}^{2t}\Big(e^{\gamma y}\overline{G_2}(x-2t)-e^{\gamma (y-t)}\overline{G_2}(x-3t)\Big)G_1(dy)\nonumber\\
& &\ \ \ \ \ \ \ \ \ \ \ \ \ \ \ \ \ \ \ \ \ \ +\int_{t}^{2t}\Big(e^{\gamma (y-t)}\overline{G_2}(x-3t)-\overline{G_2}(x-3t-(y-t))\Big)G_1(dy).
\end{eqnarray}
When $t\le y\le 2t$, we have
$$|e^{2\gamma t}\overline{G_2}(x-y)-e^{\gamma y}\overline{G_2}(x-2t)|/\overline G(x)\leq (e^{2\gamma t}/\overline{G_1}(y))+(e^{\gamma y}/\overline{G_1}(2t))$$
and
$$|e^{\gamma (y-t)}\overline{G_2}(x-3t)-\overline{G_2}(x-3t-(y-t))|/\overline G(x)\leq (e^{\gamma (y-t)}/\overline{G_1}(3t))+(1/\overline{G_1}(2t+y)).$$
And then,
$$\int_{t}^{2t}(e^{2\gamma t}/\overline{G_1}(y))+(e^{\gamma y}/\overline{G_1}(2t))G_1(dy)\leq 2e^{2\gamma t}/\overline{G_1}(2t)<\infty$$
and
$$\int_{t}^{2t}\Big((e^{\gamma (y-t)}/\overline{G_1}(3t))+(1/\overline{G_1}(2t+y)\Big)G_1(dy)\leq(e^{\gamma t}/\overline{G_1}(3t))+(1/\overline{G_1}(4t))<\infty.$$
Thus, by Fatou's lemma and (\ref{506}), we have
\begin{eqnarray}\label{5006}
&&\limsup\int_{t}^{2t}\Big(e^{2\gamma t}\overline{G_2}(x-y)-e^{\gamma y}\overline{G_2}(x-2t)\Big)G_1(dy)/\overline G(x)\nonumber\\
&\leq&\int_{t}^{2t}\limsup \Big(e^{2\gamma t}\frac{\overline{G_2}(x-y)}{\overline{G_2}(x-2t)}-e^{\gamma y}\Big)G_1(dy)/\overline G_1(2t)\leq0
\end{eqnarray}
and
\begin{eqnarray}\label{5007}
&&\limsup\int_{t}^{2t}\Big(e^{\gamma (y-t)}\overline{G_2}(x-3t)-\overline{G_2}(x-3t-(y-t))\Big)G_1(dy)/\overline G(x)\nonumber\\
&\leq&\int_{t}^{2t}\Big(\limsup\frac{e^{\gamma(y-t)}\overline{G_2}(x-3t)}{\overline{G_2}(x-3t-(y-t))}-1\Big)G_1(dy)/\overline G_1(4t)\leq0.
\end{eqnarray}
According to (\ref{502})-(\ref{5007}), we know that
\begin{eqnarray*}
\int_{t}^{2t}\Big(e^{\gamma y}\overline{G_2}(x-2t)-e^{\gamma (y-t)}\overline{G_2}(x-3t)\Big)G_1(dy)&\ge&-3\varepsilon(1+\varepsilon)e^{\gamma t}C^*(G)\overline{G}(x)-3\varepsilon \overline{G}(x).
\end{eqnarray*}
Further, by $\overline{G}(x)\approx\overline{G_2}(x-t)$, we have
\begin{eqnarray}\label{503}
\liminf\big(e^{\gamma t}\overline{G_2}(x)-\overline{G_2}(x-t)\big)/\overline{G_2}(x-t)\ge0.
\end{eqnarray}

Combined with (\ref{503}) and (\ref{506}), we know that $G_2\in\mathcal{L}(\gamma)$. \hfill$\Box$\\

Now, we prove Theorem \ref{thm40}. In Lemma \ref{lemma501}, we take $V=F,G_1=H_1,G_2=H_2=F^{*\tau}$ and $G=H$. According to Remark \ref{remark101} i), the condition (\ref{thm103}) is satisfied for the Poisson compound convolution $H_2$. Therefore, by Lemma \ref{lemma501} and $H\in\mathcal{L}(\gamma)\cap\mathcal{OS}$, we have $H_2\in\mathcal{L}(\gamma)\cap\mathcal{OS}$ and $\overline{H_2}(x)\approx\overline{H}(x)$.

Finally, by Theorem \ref{thm1}, we can get the rest of Theorem \ref{thm40}.

\section{Some local version}
\setcounter{thm}{0}\setcounter{Corol}{0}\setcounter{lemma}{0}\setcounter{pron}{0}\setcounter{equation}{0}
\setcounter{remark}{0}\setcounter{exam}{0}\setcounter{property}{0}\setcounter{defin}{0}

In this section, we give three local versions of Theorem \ref{thm40}. To this end, we first recall the concepts and notations of two local distribute classes, see, for example, Borovkov and Borovkov \cite{BB2008}.

We say that a distribution $V$ 
belongs to the distribution class $\mathcal{L}_{loc}$,
if for all $x>0$ and $0<T\le\infty$, $V(x+\Delta_T)=V(x,x+T]>0$ when $0<T<\infty$ or $V(x+\Delta_\infty)=\overline{V}(x)>0$, and for all $t>0$,
$$V(x-t+\Delta_T)\sim V(x+\Delta_T).$$

If a distribution $V$ belongs to the class $\mathcal{L}_{loc}$,
and if for all $0<T\le\infty$,
$$V^{*2}(x+\Delta_T)\sim2 V(x+\Delta_T),$$
then we say that the distribution $V$ belongs to the distribution class $\mathcal{S}_{loc}$.

Similar to the classes $\mathcal{OS}$ and $\mathcal{OL}$, we can also introduce two new distribution classes.

A distribution $V$ belongs to the class $\mathcal{OS}_{loc}$ or $\mathcal{OL}_{loc}$, if for all $0<T\le\infty$,
$$C^*_{\Delta_T}(V)=\limsup V^{*2}(x+\Delta_T)/V(x+\Delta_T)<\infty$$
or, if for all $0<T\le\infty$ and all $0<t<\infty$,
$$C^*_{\Delta_T}(V,t)=\limsup V(x+\Delta_T-t)/V(x+\Delta_T)<\infty.$$

In definitions of the above-mentioned local distribution classes, if ``for all $0<T\le\infty$" is replaced by ``for some $0<T\le\infty$", then these classes are called local long-tailed distribution class, local subexponential distribution class, generalized local long-tailed distribution class and generalized local subexponential distribution class, denoted by $\mathcal{L}_{\Delta_T}$, $\mathcal{S}_{\Delta_T}$, $\mathcal{OL}_{\Delta_T}$ and $\mathcal{OS}_{\Delta_T}$ with corresponding indicators $C^*_{\Delta_T}(V,t)$ for all $0<t<\infty$ and $C^*_{\Delta_T}(V)$, respectively. Among them, for some $0<T\le\infty$, the classes $\mathcal{L}_{\Delta_T}$ and $\mathcal{S}_{\Delta_T}$ were introduced by Asmussen et al. \cite{AFK2003}, the class $\mathcal{OS}_{\Delta_T}$ was introduced by Wang et al. \cite{WXCY2016}, and the class $\mathcal{OL}_{\Delta_T}$, as well as the classes $\mathcal{OS}_{loc}$ and $\mathcal{OL}_{loc}$, just appears in this paper.
In particular, when $T=\infty$, we get the corresponding global distribution classes $\mathcal{L}$, $\mathcal{S}$, $\mathcal{OL}$ and $\mathcal{OS}$, respectively. Compared with the definition of the class $\mathcal{S}$, however, the distribution in the class $\mathcal{S}_{\Delta_T}$ (or $\mathcal{S}_{loc}$) is required to belong to the class $\mathcal{L}_{\Delta_T}$ (or $\mathcal{L}_{loc}$), the reason of which can be found in Chen et al. \cite{CYW2013}.

\begin{thm}\label{thm6}
Let $H$ be an infinitely divisible distribution with the Laplace
transform (\ref{500}) and L$\acute{e}$vy spectral distribution $F$. For some positive constant $T_0$ and all $k\ge 1$, assume that
\begin{eqnarray}\label{thm402}
\liminf F^{*k}(x-t+\Delta_{T_0})/F^{*k}(x+\Delta_{T_0})\ge 1\ \ \text{for all}\ t>0.
\end{eqnarray}
And suppose that $H\in\mathcal{L}_{loc}\cap\mathcal{OS}_{loc}$ and $H_1(x+\Delta_{T_0})=o\big(H_2(x+\Delta_{T_0})\big)$. Then $$H_2\in\mathcal{L}_{loc}\cap\mathcal{OS}_{loc}\ and\ H_2(x+\Delta_T)\thickapprox{H}(x+\Delta_T)\ for\ all\ 0<T\le\infty.$$
And there exists a integer $n_0\ge 1$ such that $F^{*n}\in\mathcal{L}_{loc}\cap\mathcal{OS}_{loc}$ for all $n\ge n_0$.

Further, if $n_0\ge2$ and if there exits an integer $1\le l_0\le n_0-1$ such that $F^{*l_0}\in\mathcal{OS}_{loc}$, then $F^{*n}\in\mathcal{L}_{loc}\cap\mathcal{OS}_{loc}$ for all $n\ge l_0$.

In particular, if $l_0=1$, that is $F\in\mathcal{OS}_{loc}$, then $F\in\mathcal{L}_{loc}(\cap\mathcal{OS}_{loc})$.
\end{thm}

Clearly, the class $\mathcal{L}_{\Delta_T}\cap\mathcal{OS}_{\Delta_T}$ for some $0<T\le\infty$ is a larger distribution class compared with the class $\mathcal{L}_{loc}\cap\mathcal{OS}_{loc}$. And through the Esscher transform, for some $0<T<\infty$, the heavy-tailed distribution class $\mathcal{L}_{\Delta_T}\cap\mathcal{OS}_{\Delta_T}$ corresponds to the light-tailed distribution class $\mathcal{TL}_{\Delta_T}(\gamma)\cap\mathcal{OS}_{\Delta_T}$, where the class $\mathcal{TL}_{\Delta_T}(\gamma)$ is defined as follows.

For any distribution $V$ supported on $[0,\infty)$ and constant $\gamma\neq0$, if $M(V,\gamma)<\infty$, we define a distribution $V_{\gamma}$
such that
$$V_{\gamma}(x)=\int_{0-}^{x} e^{\gamma y}V(dy)\textbf{1}(x\ge0)/M(V,\gamma)$$
for all $x\in(-\infty,\infty)$, which is called the Esscher transform (or the exponential tilting) of distribution $V$.

If we consider the Esscher transform $V_{-\gamma}$ of a distribution $V$ for some $\gamma>0$, then
$$0<M(V,-\gamma)<1,\ V=(V_{-\gamma})_{\gamma}\ \text{and}\ M(V,-\gamma)M(V_{-\gamma},\gamma)=1;$$
and for all $k\ge1$,
$$(V^{*k})_{-\gamma}=(V_{-\gamma})^{*k}=V_{-\gamma}^{*k}\ \ \text{and}\ M(V^{*k},-\gamma)=\big(M(V,-\gamma)\big)^k.$$

The Esscher transformation is a key in the proofs of results of this section, because it reveals the relationship between the related heavy-tailed distribution class and light-tailed distribution class. For example, a distribution $V\in\mathcal{L}_{loc}$ if and only if $V_{-\gamma}\in\mathcal{L}(\gamma)$; and for some constant $0<T<\infty$, $V\in\mathcal{L}_{\Delta_T}$ if and only if $V_{-\gamma}$ belongs to the following distribution class
$$\mathcal{TL}_{\Delta_T}(\gamma)=\{V:\ M(V,\gamma)<\infty\ and\ V_{\gamma}\in\mathcal{L}_{\Delta_T}\}.$$
Clearly, the relationship $\mathcal{L}(\gamma)\subset\mathcal{TL}_{\Delta_T}(\gamma)$ is proper. See Definition 1.2 and Proposition 2.1 of Wang and Wang \cite{WW2011}.

Therefore, it is natural to investigate the corresponding result for the class $\mathcal{TL}_{\Delta_T}(\gamma)\cap\mathcal{OS}_{\Delta_T}$. We will find that the research method of the following result is different from the one of Theorem \ref{thm6}.

\begin{thm}\label{thm7}
Let $H$ be an infinitely divisible distribution with the Laplace transform (\ref{500}) and L$\acute{e}$vy spectral distribution $F$. For some positive constant $0<T<\infty$, assume that $H\in{\mathcal{TL}_{\Delta_T}(\gamma)}\cap\mathcal{OS}_{\Delta_T}$ for some $\gamma>0$ and $H_1(x+\Delta_{T})=o\big(H_2(x+\Delta_{T})\big)$.
In addition, suppose that for all $k\ge 1$,
\begin{eqnarray}\label{thm602}
\liminf F^{*k}_\gamma(x-t+\Delta_T)/F^{*k}_\gamma(x+\Delta_T)\ge 1\ \ for\ all\ t>0.
\end{eqnarray}
Then
$$H_2=F^{*\tau}\in{\mathcal{TL}_{\Delta_T}(\gamma)}\cap\mathcal{OS}_{\Delta_T}\ and\
H_2(x+\Delta_T)\thickapprox{H}(x+\Delta_T).$$
And there exists an integer $n_0\ge 1$ such that $F^{*n}\in\mathcal{TL}_{\Delta_T}(\gamma)\cap\mathcal{OS}_{\Delta_T}$ for all $n\ge n_0$.

Further, if $n_0\ge2$ and if there exits an integer $1\le l_0\le n_0-1$ such that $F^{*l_0}\in\mathcal{OS}_{\Delta_T}$, then $F^{*n}\in\mathcal{TL}_{\Delta_T}(\gamma)\cap\mathcal{OS}_{\Delta_T}$ for all $n\ge l_0$.

In particular, if $l_0=1$, that is $F\in\mathcal{OS}_{\Delta_T}$, then $F\in\mathcal{TL}_{\Delta_T}(\gamma)(\cap\mathcal{OS}_{\Delta_T})$.
\end{thm}



\begin{thm}\label{thm8}
Let $H$ be an infinitely divisible distribution with the Laplace transform (\ref{500}) and L$\acute{e}$vy spectral distribution $F$. For some  $0<T<\infty$, assume that condition (\ref{thm402})  with $T_0=T$ for all $k\ge 1$ is satisfied. Further, suppose that $G_1(x+\Delta_{T})=o\big(G_2(x+\Delta_{T})\big)$ and $H\in{\mathcal{L}_{\Delta_T}}\cap\mathcal{OS}_{\Delta_T}$. Then $$H_2=F^{*\tau}\in{\mathcal{L}_{\Delta_T}}\cap\mathcal{OS}_{\Delta_T}\ and\ H_2(x+\Delta_T)\thickapprox{H}(x+\Delta_T).$$

Further, if $n_0\ge2$ and if $F^{*l_0}\in\cal{OS}$ for some $1\le l_0\le n_0-1$, then $F^{*n}\in{\mathcal{L}_{\Delta_T}}\cap\mathcal{OS}_{\Delta_T}$ for all $n\ge l_0$.

In particular, if $l_0=k_0=1$, that is $F\in\mathcal{OS}_{\Delta_T}$, then $F\in\mathcal{L}_{\Delta_T}(\cap\mathcal{OS}_{\Delta_T})$.
\end{thm}

\section{The proofs of Theorem \ref{thm6}-Theorem \ref{thm8}}
\setcounter{thm}{0}\setcounter{Corol}{0}\setcounter{lemma}{0}\setcounter{pron}{0}\setcounter{equation}{0}
\setcounter{remark}{0}\setcounter{exam}{0}\setcounter{property}{0}\setcounter{defin}{0}

Here, the marks $V$, $\tau$ and $V^{*\tau}$ are the same as before.

\subsection{ Proof of Theorem \ref{thm6}}

Similar to the proof of Theorem \ref{thm40}, we first consider the closure under compound convolution roots for the distribution class $\mathcal{L}_{loc}\cap\mathcal{OS}_{loc}$. Proposition 6.1 of Xu et al. \cite{XWCY2016} note that, the class, or, more precisely, the class $(\mathcal{L}_{loc}\cap\mathcal{OS}_{loc})\setminus\mathcal{S}_{loc}$, is not closed under convolution roots.
Here, we give a positive result related compound convolution roots for the class.
\begin{thm}\label{thm2}
For any $0<\varepsilon<1$ and some constant $0<T_0<\infty$, assume that there exists an integer $n_0=n_0(V,\varepsilon,\tau,T_0)\ge1$ such that,
\begin{eqnarray}\label{thm403}
\sum\limits_{k=n_0+1}^{\infty}p_{k}V^{*(k-1)}(x+\Delta_{T_0})\leq\varepsilon V^{*\tau}(x+\Delta_{T_0})
\end{eqnarray}
for all $x\ge0$. Assume that, for all $k\ge 1$, condition (\ref{thm402}) is satisfied. If $V^{*}\in\mathcal{L}_{loc}\cap\mathcal{OS}_{loc}$, then $V^{*n}\in\mathcal{L}_{loc}\cap\mathcal{OS}_{loc}$ for all $n>n_0$.

Further, if $n_0\ge2$ and if there exists an integer $1\le l_0\le n_0-1$ such that $F^{*l_0}\in\mathcal{OS}_{loc}$, then $F^{*n}\in\mathcal{L}_{loc}\cap\mathcal{OS}_{loc}$ for all $n>l_0$.

In particular, if $l_0=1$, that is $F\in\mathcal{OS}_{loc}$, then $F\in\mathcal{L}_{loc}(\cap\mathcal{OS}_{loc})$.
\end{thm}

\proof 

In order to prove the first part of the theorem, we need the following result related to the classes $\mathcal{OS}_{\Delta_T}$ for some $0<T<\infty$ and $\mathcal{OS}_{loc}$.
\begin{lemma}\label{l401}
For any $0<\gamma,T<\infty$, $V_{-\gamma}\in\mathcal{OS}_{\Delta_T}$ if and only if $V\in\mathcal{OS}_{\Delta_T}$. Thus, for any $0<\gamma<\infty$, $V_{-\gamma}\in\mathcal{OS}_{loc}$ if and only if $V\in\mathcal{OS}_{loc}$.
\end{lemma}
\proof From (2.4) of Wang and Wang \cite{WW2011} that, for any $0<\gamma,T<\infty$,
\begin{eqnarray}\label{4070}
&&V_{-\gamma}(x+\Delta_{T})=M(V_{-\gamma},\gamma)e^{-\gamma x}V(x+\Delta_{T})\Big(1-\gamma\int_0^{T}\frac{V(x+y,x+T]}{e^{\gamma y}
V(x+\Delta_{T})}dy\Big),
\end{eqnarray}
where $V(x+\Delta_{T})>0$ for all $x\ge0$, we obtain the following inequality,
\begin{eqnarray}\label{40700}
e^{-\gamma T}V(x+\Delta_{T})\le e^{\gamma x}V_{-\gamma}(x+\Delta_{T})/M(V_{-\gamma},\gamma)\le V(x+\Delta_{T}).
\end{eqnarray}
Thus, if $V\in\mathcal{OS}_{\Delta_T}$, then by (\ref{40700}) and Radon-Nikodym Theorem, we have
\begin{eqnarray*}
&&V^{*2}_{-\gamma}(x+\Delta_{T})=\int_0^x V_{-\gamma}(x-y+\Delta_{T})V_{-\gamma}(dy)
+\int_x^{x+T}V_{-\gamma}(0,x-y+T]V_{-\gamma}(dy)\nonumber\\
&\le&e^{-\gamma x}M(V_{-\gamma},\gamma)\int_0^x V(x-y+\Delta_{T})V(dy)/M(V,-\gamma)+V_{-\gamma}(x+\Delta_T)\nonumber\\
&\le&e^{-\gamma x}M(V_{-\gamma},\gamma)V^{*2}(x+\Delta_{T})/M(V,-\gamma)+V_{-\gamma}(x+\Delta_T)\nonumber\\
&\le&2C^*_{\Delta_T}(V)M(V_{-\gamma},\gamma)e^{-\gamma x}V(x+\Delta_T)/M(V,-\gamma)+V_{-\gamma}(x+\Delta_T)\nonumber\\
&\le&\big(2C^*_{\Delta_T}(V)M(V_{-\gamma},\gamma)+1\big)e^{\gamma T}V_{-\gamma}(x+\Delta_T)/M(V,-\gamma),
\end{eqnarray*}
that is $V_{-\gamma}\in\mathcal{OS}_{\Delta_T}$. Conversely, if $V_{-\gamma}\in\mathcal{OS}_{\Delta_T}$, we can also get $V\in\mathcal{OS}_{\Delta_T}$ by the same way. Therefore, $V_{-\gamma}\in\mathcal{OS}_{\Delta_T}$ if and only if $V\in\mathcal{OS}_{\Delta_T}$ for some $0<T<\infty$. According to the arbitrariness of $T$, we can also  prove that $V_{-\gamma}\in\mathcal{OS}_{loc}$ if and only if $V\in\mathcal{OS}_{loc}$.
\hfill$\Box$\\

Now, we continue to prove Theorem \ref{thm2}. In Lemma \ref{l401}, we replace $V$ with $V^{*\tau}$. Then by $V^{*\tau}\in\mathcal{L}_{loc}\cap\mathcal{OS}_{loc}$ and Proposition 2.2 of Wang et al. \cite{WXCY2016}, we know that
$$(V^{*\tau})_{-\gamma}\in\mathcal{L}(\gamma)\cap\mathcal{OS}_{loc}\subset\mathcal{L}(\gamma)\cap\mathcal{OS}.$$
In addition,
$$M(V^{*\tau},-\gamma)=\sum_{k=0}^\infty p_kM^k(V,-\gamma)=E(M(F,-\gamma))^\tau<1$$
and for all $x\ge0$,
\begin{eqnarray}\label{thm404}
\overline{(V^{*\tau})_{-\gamma}}(x)&=&\sum_{k=1}^\infty p_kM^k(V,-\gamma) \overline{V_{-\gamma}^{*k}}(x)/M(V^{*\tau},-\gamma)\nonumber\\
&=&\sum_{k=1}^\infty q_k\overline{V_{-\gamma}^{*k}}(x)=\overline{(V_{-\gamma})^{*\sigma}}(x).
\end{eqnarray}

For any $0<\varepsilon<1$, there is a number $0<\varepsilon_0\le e^{-T_0}$ such that $\varepsilon=\varepsilon_0e^{\gamma T_0}$. By (\ref{40700}), (\ref{thm404}) and (\ref{thm403}) with $\varepsilon_0$ and the corresponding $n_0=n_0(V,\varepsilon_0,\tau,T_0)$,  we have
\begin{eqnarray}\label{thm408}
&&\sum_{k=n_0+1}^\infty q_k\overline{V^{*(k-1)}_{-\gamma}}(x)=\sum_{k=n_0+1}^\infty q_k \sum_{m=0}^\infty¡¡V_{-\gamma}^{*(k-1)}(x+mT_0+\Delta_{T_0})\nonumber\\
&\le&\sum_{k=n_0+1}^\infty p_k(M(V,-\gamma))^{k-1} (M(V_{-\gamma},\gamma))^{k-1} \sum_{m=0}^\infty¡¡ e^{-\gamma(x+mT_0)}V^{*(k-1)}(x+mT_0+\Delta_{T_0})/M(V^{*\tau},-\gamma)\nonumber\\
&=&\sum_{k=n_0+1}^\infty p_k \sum_{m=0}^\infty e^{-\gamma(x+mT_0)}V^{*(k-1)}(x+mT_0+\Delta_{T_0})/M(V^{*\tau},-\gamma)\nonumber\\
&=&\sum_{m=0}^\infty e^{-\gamma(x+mT_0)}\sum_{k=n_0+1}^\infty p_k V^{*(k-1)}(x+mT_0+\Delta_{T_0})/M(V^{*\tau},-\gamma)\nonumber\\
&<&\varepsilon_0\sum_{m=0}^\infty e^{-\gamma(x+mT_0)}F^{*\tau}(x+mT_0+\Delta_{T_0})/M(F^{*\tau},-\gamma)\nonumber\\
&\le&\varepsilon_0 e^{\gamma T_0}\sum_{m=0}^\infty(V^{*\tau})_{-\gamma}(x+mT_0+\Delta_{T_0})\nonumber\\
&=&\varepsilon\overline{(V^{*\tau})_{-\gamma}}(x)=\varepsilon\overline{(V_{-\gamma})^{*\sigma}}(x)
\end{eqnarray}
for all $x\ge0$. Further, by (\ref{thm402}) and Proposition 2.1 of Wang and Wang \cite{WW2011}, for all $k\ge k_0$,
\begin{eqnarray}\label{thm410}
\liminf \overline{V^{*k}_{-\gamma}}(x-t)/\overline{V^{*k}_{-\gamma}}(x)\ge e^{\gamma t}\ for\ all\ t>0.
\end{eqnarray}

Since $(V^{*\tau})_{-\gamma}\in\mathcal{L}(\gamma)\cap\mathcal{OS}$, by (\ref{thm408}), (\ref{thm410}) and Theorem \ref{thm1}, we have $V^{*n}_{-\gamma}\in\mathcal{L}(\gamma)\cap\mathcal{OS}$ for all $n\ge n_0$. We note that, when $V^{*n}_{-\gamma}\in\mathcal{L}(\gamma)$, for any $0<T<\infty$,
\begin{eqnarray}\label{thm411}
V^{*n}_{-\gamma}(x+\Delta_T)\sim(1-e^{-\gamma T})\overline{V^{*n}_{-\gamma}}(x).
\end{eqnarray}
Thus, $V^{*n}_{-\gamma}\in\mathcal{L}(\gamma)\cap\mathcal{OS}=\mathcal{L}(\gamma)\cap\mathcal{OS}_{loc}$, that is $V^{*n}\in\mathcal{L}_{loc}\cap\mathcal{OS}_{loc}$ for all $n\ge n_0$.

In the same way, we can prove the rest of the theorem. \hfill$\Box$\\

In order to prove Theorem \ref{thm6}, we also need the following lemma. Recall that $G_1$ is a distribution, $G_2=V^{*\tau}$ and $G=G_1*G_2$. By Esscher transform, Lemma \ref{lemma501} and Lemma \ref{l401}, we can prove the lemma.

\begin{lemma}\label{lemma502}
For some constant $T_0>0$ and any $0<\varepsilon<1$, assume that conditions (\ref{thm403}) and (\ref{thm402}) for all $k\ge 1$ and some $\gamma\ge0$  are satisfied. Further, suppose that $G_1(x+\Delta_{T_0})=o\big(G_2(x+\Delta_{T_0})\big)$ and $G\in{\mathcal{L}_{loc}}\cap\mathcal{OS}_{loc}$. Then $G_2\in{\mathcal{L}_{loc}}\cap\mathcal{OS}_{loc}$ and
$G_2(x+\Delta_{T_0})\thickapprox G(x+\Delta_{T_0})$.
\end{lemma}

Clearly, condition (\ref{thm403}) is satisfied for any $0<T\le\infty$, if $p_n=e^{-\lambda}\lambda^n/n!$ for all nonnegative integers $n$. Therefore, by Lemma \ref{lemma502} and Theorem \ref{thm2}, we can prove Theorem \ref{thm6}.

\subsection{ Proof of Theorem \ref{thm7}}

According to Lemma \ref{lemma502}, we only need to prove the following result.
\begin{thm}\label{thm601}
Assume that $M(V,\gamma)<\infty$ for some $0<\gamma<\infty$ and for any $0<\varepsilon<1$ and some $0<T<\infty$, there exists an integer $n_0=n_0(F,\tau,\varepsilon,T)\ge 1$ such that (\ref{thm403}) with $T_0=T$ is holds for all $x\ge0$.
In addition, suppose that condition (\ref{thm602}) for all $k\ge 1$ is satisfied.
If $V^{*\tau}\in{\mathcal{TL}_{\Delta_T}(\gamma)}\cap\mathcal{OS}_{\Delta_T}$ for some $0<\gamma<\infty$, then 
$V^{*n}\in{\mathcal{TL}_{\Delta_T}(\gamma)}\cap\mathcal{OS}_{\Delta_T}$ for all $n\ge n_0$.

Further, if $n_0\ge2$ and if there exits an integer $1\le l_0\le n_0-1$ such that $V^{*l_0}\in\mathcal{OS}_{\Delta_T}$, then $V^{*n}\in\mathcal{TL}_{\Delta_T}\cap\mathcal{OS}_{\Delta_T}$ for all $n\ge l_0$.

In particular, if $l_0=1$, that is $V\in\mathcal{OS}_{\Delta_T}$, then $V\in\mathcal{TL}_{\Delta_T}(\cap\mathcal{OS}_{\Delta_T})$.
\end{thm}
\begin{remark}\label{remark602}
For any $0<\gamma,T<\infty$, from (\ref{thm411}), the condition in Theorem \ref{thm601} that $F^{*\tau}\in{\mathcal{TL}_{\Delta_T}(\gamma)}\cap\mathcal{OS}_{\Delta_T}$ is substantially weaker than the corresponding condition in Theorem \ref{thm1} and Theorem \ref{thm2} that $V^{*\tau}\in{\mathcal{L}}(\gamma)\cap\mathcal{OS}={\mathcal{L}}(\gamma)\cap\mathcal{OS}_{loc}$.
\end{remark}

\proof 
In order to prove the theorem, we first give the local version similar to the half of Lemma \ref{lemma1}.

\begin{lemma}\label{lemma2}
If $V^{*\tau}\in\mathcal{OS}_{\Delta_T}$ for some $0<T<\infty$, then the following proposition ii) can deduce the proposition i):

i) For any $0<\varepsilon<1$, there is an integer $n_0=n_0(F,\tau,\varepsilon,T)\ge1$ such that
\begin{eqnarray}\label{l501}
\sum\limits_{k=n_0+1}^{\infty}p_kV^{*k}(x+\Delta_T)\leq\varepsilon V^{*\tau}(x+\Delta_T)\ for\ all\ x\ge0.
\end{eqnarray}

ii) For any $0<\varepsilon<1$, there is an integer $n_0=n_0(V,\tau,\varepsilon,T)\ge1$ such that
\begin{eqnarray}\label{l502}
\sum\limits_{k=n_0+1}^{\infty}p_kV^{*(k-1)}(x+\Delta_T)\leq\varepsilon V^{*\tau}(x+\Delta_T)\ for\ all\ x\ge0.
\end{eqnarray}
\end{lemma}
\proof  In order to prove ii) $\Longrightarrow$ i), we denote
$$D^*(V^{*\tau},T)=\sup_{x\ge0}(V^{*\tau})^{*2}(x+\Delta_T)/V^{*\tau}(x+\Delta_T).$$
Clearly, $0< D^*(V^{*\tau},T)<\infty$. For any $0<\varepsilon<1$, there is a positive number $\varepsilon_1$ such that $$0<\varepsilon_1<p_1/\big(1+D^*(V^{*\tau},T)\big)<1,$$
and $\varepsilon=\varepsilon_1\big(1+D^*(V^{*\tau},T)\big)/p_1.$
For above $\varepsilon_1$, by proposition ii), there is an integer $n_0=n_0(V,\tau,\varepsilon_1,T)\ge1$ such that (\ref{l502}) holds.
Further, we assume that
\begin{eqnarray}\label{l103}
a_{n_0}=\sum\limits_{k=n_0+1}^{\infty}p_k<\varepsilon_1.
\end{eqnarray}
Then by (\ref{l501}) and (\ref{l103}), for all $x\ge0$, we have
\begin{eqnarray*}
&&\sum\limits_{k=n_0+1}^{\infty}p_{k}V^{*k}(x+\Delta_T)=V*\sum\limits_{n_0+1\le k<\infty}p_{k}V^{*(k-1)}(x+\Delta_T)\nonumber\\ &\le&a_{n_0}\int_{0-}^xG_{n_0}(x-y+\Delta_T)V(dy)+a_{n_0}V(x+\Delta_T)\nonumber\\ &<&\varepsilon_1\Big(\int_{0-}^xV^{*\tau}(x-y+\Delta_T)V^{*\tau}(dy)+V^{*\tau}(x+\Delta_T)\Big)/p_1\nonumber\\ &\le&\varepsilon_1\Big(V^{*2\tau}(x+\Delta_T)+V^{*\tau}(x+\Delta_T)\Big)/p_1\nonumber\\
&\le&\varepsilon_1 \big(1+D^*(V^{*\tau},T)\big) V^{*\tau}(x+\Delta_T)/p_1 =\varepsilon V^{*\tau}(x+\Delta_T),
\end{eqnarray*}
that is we get proposition i). \hfill$\Box$\\

In the following, we continue to prove the theorem. First, we prove that $V^{*n}\in\mathcal{OS}_{\Delta_T}$ for all $n\ge n_0.$

Because $V^{*\tau}\in{\mathcal{TL}_{\Delta_T}(\gamma)}$, $M(V^{*\tau},\gamma)<\infty$ which implies $M(V^{*n},\gamma)<\infty$ for all $n\ge1$. From (2.4) of Wang and Wang \cite{WW2011}, it holds that
\begin{eqnarray}\label{601}
V^{*\theta}(x+\Delta)=M(V^{*\theta},\gamma)e^{-\gamma x}(V^{*\theta})_\gamma(x+\Delta_{T})
\Big(1-\gamma\int_0^{T}\frac{(V^{*\theta})_\gamma(x+y,x+T]}
{e^{\gamma y}(V^{*\theta})_\gamma(x+\Delta_{T})}dy\Big),
\end{eqnarray}
where $\theta=k$ for all $k\ge1$ or, $\theta=\tau$. In particular, when $\theta=\tau$,
\begin{eqnarray}\label{602}
(V^{*\tau})_\gamma(x+\Delta_{T})&=&\sum_{k=1}^\infty p_kM^k(V,\gamma) V_{\gamma}^{*k}(x+\Delta_{T})/M(V^{*\tau},\gamma)
=(V_{\gamma})^{*\sigma}(x+\Delta_{T})
\end{eqnarray}
for all $x\ge0$, where $\sigma$ is a random variable such that
$$q_k=P(\sigma=k)=p_kM^k(V,\gamma)/M(V^{*\tau},\gamma)$$
for all $k\ge0$. Thus,
\begin{eqnarray}\label{603}
e^{-\gamma T}(V^{*\theta})_\gamma(x+\Delta_{T})
\le e^{\gamma x}V^{*\theta}(x+\Delta)/M(V^{*\theta},\gamma)\le(V^{*\theta})_\gamma(x+\Delta_{T}).
\end{eqnarray}
For any function $h$ such that $0<h(x)\uparrow\infty$ and $h(x)/x\rightarrow0$, by (\ref{603}) and Radon-Nikodym theorem, we have
\begin{eqnarray*}
&&\int_{h(x)}^{x-h(x)}(V^{*\theta})_\gamma(x-y+\Delta_{T})(V^{*\theta})_\gamma(dy)
=O\Big(e^{\gamma x}\int_{h(x)}^{x-h(x)}V^{*\theta}(x-y+\Delta_{T})V^{*\theta}(dy)\Big).
\end{eqnarray*}
Here the notation $f(x)=O(g(x))$ for positive-valued
functions $f$ and $g$ means $\limsup f(x)/g(x)\\<\infty$. Combined with the asymptotic inequality and Lemma \ref{l401}, when $V^{*\tau}\in{\mathcal{TL}_{\Delta_T}(\gamma)}\cap\mathcal{OS}_{\Delta_T}$,
$(V^{*\tau})_\gamma\in\mathcal{L}_{\Delta_T}\cap\mathcal{OS}_{\Delta_T}$.

According (\ref{thm403}), (\ref{603}) and Lemma \ref{lemma2}, for any $0<\varepsilon<1$, there is an integer $n_0=n_0(V,\tau,\varepsilon,\gamma,T)$ such that for\ all\ $x\ge0$
\begin{eqnarray}\label{604}
\sum\limits_{k=1}^{n_0}q_{k}V^{*k}_\gamma(x+\Delta_T)&\ge& e^{\gamma x}\sum\limits_{k=1}^{n_0}p_{k}V^{*k}(x+\Delta_T)/M(V^{*\tau},\gamma)\nonumber\\
&\ge&(1-\varepsilon)e^{\gamma x} V^{*\tau}(x+\Delta_T)/M(V^{*\tau},\gamma)\nonumber\\
&\ge&(1-\varepsilon)e^{-\gamma T} (V^{*\tau})_\gamma(x+\Delta_T).
\end{eqnarray}
Further, for any $n\ge k_0$, used Fatou's lemma, Radon-Nikodym theorem, (\ref{603}) and (\ref{thm602}), respectively, we have
\begin{eqnarray}\label{606}
\liminf \frac{V^{*(n+1)}_\gamma(x+\Delta_T)}{V^{*n}_\gamma(x+\Delta_T)}
&\ge&\int_0^{\infty}\liminf\frac{V^{*n}_\gamma(x-y+\Delta_T)}{V^{*n}_\gamma(x+\Delta_T)}\textbf{1}(y\le x)F_\gamma(dy)
=1.
\end{eqnarray}

Combined with (\ref{604}), (\ref{606}) and $(V^{*\tau})_\gamma\in\mathcal{OS}_{\Delta_T}$,
we know that $V^{*n_0}_\gamma\in\mathcal{OS}_{\Delta_T}$ and $V^{*n_0}_\gamma(x+\Delta_T)\approx (V^{*\tau})_\gamma(x+\Delta_T)$, where the notation $f(x)\thickapprox g(x)$ for positive-valued
functions $f$ and $g$ means $f(x)=O\big(g(x)\big)$ and $g(x)=O\big(f(x)\big)$.
Use again Lemma \ref{l401}, (\ref{40700}) and (\ref{thm602}), for all $n\ge n_0$, we have
$$V^{*n}\in\mathcal{OS}_{\Delta_T}\ and\ V^{*n}(x+\Delta_T)\approx V^{*\tau}(x+\Delta_T).$$

Next, we prove that $V^{*n}\in{\mathcal{TL}_{\Delta_T}(\gamma)}$ for all $n\ge n_0$.

According to Lemma \ref{lemma2}, (\ref{thm403}), (\ref{602}) and (\ref{603}), for any $0<\varepsilon<1$ and any fixed $n\ge n_0$, there exists an integer $m_0=m_0(F,\tau,\varepsilon,T)\ge n$ such that
\begin{eqnarray}\label{t103}
\sum\limits_{k=m_0+1}^{\infty}q_kV^{*k}_\gamma(x+\Delta_T)\leq\varepsilon(V^{*\tau})_\gamma(x+\Delta_T)\ \text{for all}\ x\ge0.
\end{eqnarray}

Further, by $(V^{*\tau})_\gamma\in\mathcal{L}_{\Delta_T}\cap\mathcal{OS}_{\Delta_T}$ and (\ref{thm602}), for the above $\varepsilon$ and any $t>0$, there is a constant $x_0=x_0(V,\tau,\varepsilon,t)$ such that, for all $x>x_0$,
\begin{eqnarray*}
&&\varepsilon (V^{*\tau})_\gamma(x+\Delta_T)\ge (V^{*\tau})_\gamma(x-t+\Delta_T)-(V^{*\tau})_\gamma(x+\Delta_T)\nonumber\\
&=&\Big(\sum\limits_{1\le k\neq n\le m_0}+\sum\limits_{k=n}+\sum\limits_{k\ge m_0+1}\Big)q_k\big(V^{*k}_\gamma(x-t+\Delta_T)-V^{*k}_\gamma(x+\Delta_T)\big)\nonumber\\
&\geq&-\varepsilon \sum\limits_{k_0\le k\le m_0}q_kV^{*k}_\gamma(x+\Delta_T)+q_{n}\big(V^{*n}_\gamma(x-t+\Delta_T)-V^{*n}_\gamma(x+\Delta_T)\big)-\varepsilon(V^{*\tau})_\gamma(x+\Delta_T),
\end{eqnarray*}
which implies that, for all $x>x_{0}$,
\begin{eqnarray*}
V^{*n}_\gamma(x-t+\Delta_T)&\le& V^{*n}_\gamma(x+\Delta_T)+3\varepsilon (V^{*\tau})_\gamma(x+\Delta_T)/q_{n}.
\end{eqnarray*}
Hence by $(V^{*\tau})_\gamma(x+\Delta_T)\approx V^{*n}_\gamma(x+\Delta_T)$, $V^{*\tau}\in{\mathcal{TL}_{\Delta_T}(\gamma)}$ and arbitrariness of $\varepsilon$, we can get
\begin{eqnarray}\label{t104}
\limsup V^{*{n}}_\gamma(x-t+\Delta_T)/V^{*{n}}_\gamma(x+\Delta_T)\le 1.
\end{eqnarray}
Combining (\ref{t104}) and (\ref{thm602}), $V^{*n}_\gamma\in{\cal{L}}_{\Delta_T}$. Therefore, $V^{*n}\in{\mathcal{TL}_{\Delta_T}(\gamma)}$, for all $n\ge n_0$.

Finally, the remainder of the theorem can be similarly proved.
\hfill$\Box$\\

\subsection{ Proof of Theorem \ref{thm8}}

Along the same way of Theorem \ref{thm601} without the Esscher transform, we can get the following conclusion and omit its proof.
\begin{thm}\label{thm6020}
Assume that for any $0<\varepsilon<1$ and some $0<T<\infty$, there exists an integer $n_0=n_0(V,\tau,\varepsilon,T)\ge 1$ such that the condition (\ref{thm403}) with $T_0=T$ is satisfied. In addition, suppose that (\ref{thm402}) with $T_0=T$ for all $k\ge k_0$ is satisfied. If $V^{*\tau}\in{\mathcal{L}_{\Delta_T}}\cap\mathcal{OS}_{\Delta_T}$, then
$V^{*n}\in{\mathcal{L}_{\Delta_T}}\cap\mathcal{OS}_{\Delta_T}$ for all $n\ge n_0$.

Further, if $l_0\ge2$ and if $V^{*l_0}\in\cal{OS}$ for some $1\le l_0\le n_0-1$, then $V^{*n}\in{\mathcal{L}_{\Delta_T}}\cap\mathcal{OS}_{\Delta_T}$ for all $n\ge l_0$.

In particular, if $l_0=1$, that is $V\in\mathcal{OS}_{\Delta_T}$, then $V\in\mathcal{L}_{\Delta_T}(\cap\mathcal{OS}_{\Delta_T})$.
\end{thm}

Based on Theorem \ref{thm6020}, we can get result of Theorem \ref{thm8} without its details of proof.

\section{Discussion}
\setcounter{thm}{0}\setcounter{Corol}{0}\setcounter{lemma}{0}\setcounter{pron}{0}\setcounter{equation}{0}
\setcounter{remark}{0}\setcounter{exam}{0}\setcounter{property}{0}\setcounter{defin}{0}


In this section, we follow the notations of Theorem \ref{thm40} and Theorem \ref{thm1}, respectively.

\subsection{On the condition (\ref{thm102})}

In this subsection, we provide some distributions satisfying the condition (\ref{thm102}) for all $k\ge 1$.  %

\begin{pron}\label{p0}
i) For $i=1,2$, let $F_i$ be a distribution satisfying
\begin{eqnarray}\label{l201}
\liminf\overline{F_i}(x-t)/\overline{F_i}(x)\ge e^{\gamma t}\ for\ any\ t>0.
\end{eqnarray}
Assume that $F_2\in\mathcal{OL}$ and
\begin{eqnarray}\label{l202}
\lim\overline{F_1}(x)C^*(F_2,x)=0.
\end{eqnarray}
Then
\begin{eqnarray}\label{l203}
\liminf\overline{F_1*F_2}(x-t)/\overline{F_1*F_2}(x)\ge e^{\gamma t}\ for\ all\ t>0.
\end{eqnarray}

ii) Further, if $F_1=F, F_2=F^{*k_0}\in\mathcal{OL}$ for some integer $k_0\ge1$, then (\ref{thm102}) holds for all $k\ge k_0$.
In particular, if $k_0=1$, then (\ref{thm102}) holds for all $k\ge1$.
\end{pron}

\begin{remark}\label{remark601}

The condition (\ref{l202}) is necessary in certain sense. In Proposition 5.4 of Xu et al. \cite{XWCY2016}, there is a distribution $F$ such that $F\in\mathcal{OL}\setminus\mathcal{L}(\gamma)$ for each $\gamma>0$ and $F^{*k}\in\mathcal{L}(\gamma)$ for all $k\ge2$, thus $k_0=1$. However, for $k=1$ the condition (\ref{thm102}) is not satisfied, and for $k\ge2$ it is satisfied. Clearly, the condition (\ref{l202}) is not holds, otherwise, (\ref{thm102}) holds for all $k\ge1$ by Proposition \ref{p0} ii).
\end{remark}
\proof i) To prove (\ref{l203}), we only need to prove its equivalent proposition:
\begin{eqnarray}\label{l2050}
\limsup(e^{\gamma t}\overline{F_1*F_2}(x)-\overline{F_1*F_2}(x-t))/\overline{F_1*F_2}(x-t)\le0.
\end{eqnarray}
From (\ref{l201}), we know that, for any $0<\varepsilon<1$ and $i=1,2$, there exists a constant $x_0=x_0(\varepsilon,\gamma,F_i,i=1,2)>0$ such that for all $x>x_0$
\begin{eqnarray}\label{l205}
\overline{F_i}(x-t)-e^{\gamma t}\overline{F_i}(x)\geq-\varepsilon\overline{F_i}(x).
\end{eqnarray}
For any $t>0$, when $x>t+2x_0$, by (\ref{l205}), we have
\begin{eqnarray*}
&&e^{\gamma t}\overline{F_1*F_2}(x)-\overline{F_1*F_2}(x-t)\nonumber\\
&=&e^{\gamma t}\int_{x-t-x_0}^{x-x_0}\overline{F_1}(x-y)F_2(dy)-\overline{F_1}(x_0)(\overline{F_2}(x-t-x_0)-e^{\gamma t}\overline{F_2}(x-x_0))\nonumber\\
& &-\int_{0}^{x-t-x_0}\Big(\overline{F_1}(x-t-y)-e^{\gamma t}\overline{F_1}(x-y)\Big)F_2(dy)-\int_{0}^{x_0}\Big(\overline{F_2}(x-t-y)-e^{\gamma t}\overline{F_2}(x-y)\Big)F_1(dy)
\nonumber\\
&\leq&e^{\gamma t}\overline{F_1}(x_0)\overline{F_2}(x-t-x_0)+\varepsilon\overline{F_1}(x_0) \overline{F_2}(x-x_0)\nonumber\\
& &+\varepsilon\int_{0}^{x-t-x_0}\overline{F_1}(x-y)F_2(dy)+\varepsilon\int_{0}^{x_0} \overline{F_2}(x-y)F_1(dy)\nonumber\\
&\leq&e^{\gamma t}\overline{F_1}(x_0)\overline{F_2}(x-t-x_0)+\varepsilon\overline{F_1*F_2}(x).
\end{eqnarray*}
Then
$$\limsup\big(e^{\gamma t}\overline{F_1*F_2}(x)-\overline{F_1*F_2}(x-t)\big)/\overline{F_1*F_2}(x-t)\le e^{\gamma t}\overline{F_1}(x_0)C^*(F_2,x_0)+\varepsilon.$$
Therefore, by (\ref{l202}) and arbitrariness of $\varepsilon$, (\ref{l2050}) holds.

ii) We first prove the following fact: if $F\in\mathcal{OL}$, then for all $k\ge1$, $F^{*k}\in\mathcal{OL}$ and
\begin{eqnarray}\label{l206}
C^*(F^{*k},t)\le C^*(F,t)\ for\ all\ t\ge0.
\end{eqnarray}
We use mathematical induction to prove the result. Clearly, (\ref{l206}) holds for $k=1$. Assume that $F^{*k}\in\mathcal{OL}$ for some $k\ge1$, then for any $0<\varepsilon<1$ and $t>0$, there is a constant $x_0=x_0(F^{*k},\varepsilon,t)$ such that when $x\ge x_0$,
$$\overline{F^{*k}}(x-t)\le (1+\varepsilon)C^*(F^{*k},t)\overline{F^{*k}}(x)\le (1+\varepsilon)C^*(F,t)\overline{F^{*k}}(x)<\infty.$$
Further, according to the induction hypothesis, for any $0<\varepsilon<1$ and $t>0$, we have
\begin{eqnarray*}\label{l207}
&&\overline{F^{*(k+1)}}(x-t)=\Big(\int_{0-}^{x-t-x_0}+\int_{x-t-x_0}^{x-t}\Big)\overline{F^{*k}}(x-t-y)F(dy)+\overline{F}(x-t)\nonumber\\
&\le&(1+\varepsilon)C^*(F,t)\Big(\int_{0-}^{x-x_0}\overline{F^{*k}}(x-y)F(dy)+\int_{0-}^{x_0}\overline{F}(x-y)F^{*k}(dy)
+\overline{F}(x-x_0)\overline{F^{*k}}(x_0)\Big)\nonumber\\
&\le&(1+\varepsilon)C^*(F,t)\overline{F^{*(k+1)}}(x).
\end{eqnarray*}
Thus, $F^{*(k+1)}\in\mathcal{OL}$ and (\ref{l206}) holds for $k+1$ by the arbitrariness of $\varepsilon$.

Next, for any $m\ge k_0+1$, we take $F_1=F$ and $F_2=F^{*(m-1)}$, then by (\ref{l206}) and (\ref{l202}), we have
$$\overline{F_1}(x)C^*(F_2,x)=\overline{F}(x)C^*(F^{*(m-1)},x)\le\overline{F}(x)C^*(F,x)\to0.$$
Therefore, by the conclusion in 1), (\ref{thm102}) holds for $k=m$.
\hfill$\Box$\\

Now, we will introduce a kind of distribution with more specific representation, which satisfies (\ref{thm102}).
For some constant $\gamma>0$ and arbitrary distribution $F_0$, we define the distribution $F$ in the form
\begin{eqnarray}\label{0301}
\overline{F}(x)=\textbf{1}(x<0)+e^{-\gamma x}\overline{F_0}(x)\textbf{1}(x\ge0),\ x\in(-\infty,\infty).
\end{eqnarray}
Clearly, $F$ is light-tailed and (\ref{thm102}) holds for $k=k_0=1$.

\begin{pron}\label{p1}
In (\ref{0301}), if $F_0$ is a heavy-tailed distribution such that $F_0\in\mathcal{OL}$ and
\begin{eqnarray}\label{l301}
\lim\overline{F_0}(x)C^*(F_0,x)=0,
\end{eqnarray}
then $F^{*k}\in\mathcal{OL}$ and (\ref{thm102}) holds for all $k\ge 1$.
\end{pron}
\proof~~Clearly, (\ref{thm102}) holds for $k=1$ and
\begin{eqnarray*}
C^*(F,t)&=&\limsup e^{\gamma t}\overline{F_0}(x-t)/\overline{F_0}(x)=e^{\gamma t}C^*(F_0,t).
\end{eqnarray*}
Thus, $F\in\mathcal{OL}$, and by (\ref{l301}),
$$\overline{F}(x)C^*(F,x)=\overline{F_0}(x)C^*(F_0,x)\to0.$$
Therefore, by Proposition \ref{p0}, (\ref{thm102}) holds for all $k\ge1$.\hfill$\Box$\\

Finally, we give a distribution $F$ satisfying the conditions (\ref{l202}) and (\ref{thm102}).
\begin{exam}\label{exam-01}
Let $\alpha\in(3/2,(\sqrt{5}+1)/2)$ and $r=(\alpha+1)/\alpha$ be constants.
Assume $a>1$ is enough large such  that $a^{r} > 8a$.
We define a distribution $F_0$ supported on $[0,\infty)$ such that
\begin{eqnarray}\label{exam01}
\overline{F_0}(x)&=&\textbf{\emph{1}}(x<
a_0)+C\sum\limits_{n=0}^{\infty}\Big(\big(\sum\limits_{i=n}^{\infty}
a_i^{-\alpha}-a_n^{-\alpha-1}(x-a_n)\big)
\textbf{\emph{1}}\big(x\in[a_n,2a_n)\big)\nonumber\\
& &
+\sum\limits_{i=n+1}^{\infty}a_i^{-\alpha}\textbf{\emph{1}}\big(x\in[2a_n,
a_{n+1})\big)\Big),
\end{eqnarray}
where $C$ is a regularization constant and $a_n=a^{r^n}$ for all nonnegative integers, see Proposition 4.3 of Xu et al. \cite{XWCY2016}.

Let $F$ be a distribution defined by (\ref{0301}). Then, for any $t>0$ and all enough large integer $n$ such that $2a_n+t<a_{n+1}$,
when $x\in[a_n,a_n+t)$,
\begin{eqnarray*}
e^{\gamma t}<\overline {F}(x-t)/\overline {F}(x)&=&e^{\gamma t}\overline {F_0}(x-t)/\overline {F_0}(x)
\le e^{\gamma t}\overline {F_0}(a_n)/\overline {F_0}(a_n+t)\to e^{\gamma t};
\end{eqnarray*}
when $x\in[a_n+t,2a_n)$,
\begin{eqnarray*}
e^{\gamma t}<\overline {F}(x-t)/\overline {F}(x)&=&e^{\gamma t}\big(\overline {F_0}(x)-a_n^{-\alpha-1}t\big)/\overline {F_0}(x)
\leq e^{\gamma t}\big(\overline {F_0}(a_n)-a_n^{-\alpha-1}t\big)/\overline {F_0}(a_n)\to e^{\gamma t};
\end{eqnarray*}
when $x\in[2a_n,2a_n+t)$, from $r=(\alpha+1)/\alpha$,
$$\overline {F}(x-t)/\overline {F}(x)\le e^{\gamma t}\overline {F_0}(2a_n-t)/\overline {F_0}(2a_n)\to e^{\gamma t}(1+t);$$
and when $x\in[2a_n+t,a_{n+1})$,
$$\overline {F}(x-t)/\overline {F}(x)=e^{\gamma t}.$$
This fact implies $C^*(F,t)=(1+t)e^{\gamma t}$, $F\notin\mathcal{L}(\gamma)$ and $F^{*k}\in\mathcal{OL}$ for all $k\ge1$.
Further,  by $\int_0^\infty \overline{F}(y)dy<\infty$, we have
\begin{eqnarray*}
&&\lim\overline{F}(x)C^*(F,x)=\lim\overline{F_0}(x)(1+x)=0.
\end{eqnarray*}
Therefore, by Proposition \ref{p1}, (\ref{thm102}) holds for all $k\ge1$.\hfill$\Box$\\
\end{exam}

\subsection{On the condition (\ref{thm103})}

In this subsection, we give a more general Kesten inequality, by which, we can implies (\ref{thm103}) under certain conditions. To this end, write $A_n=\sup_{x\geq0}\frac{\overline{V^{*n}}(x)}{\overline{G}(x)}$ for all $n\ge1$ and
$$a=M(V,\gamma)+A_1\big(C^*(G)-2M(G,\gamma)\big)=M(V,\gamma)+b.$$
\begin{pron}\label{p2}
Let $V$ and $G$ be two distributions such that $G\in\mathcal{L}(\gamma)\cap\mathcal{OS}$ for some $0\le\gamma<\infty$, $M(G,\gamma)<\infty$ and $\overline{V}(x)=O(\overline{G}(x))$. Then for any constant $M$ satisfying
\begin{eqnarray}\label{p20110}
a<M<1+a
\end{eqnarray}
and any constant $\varepsilon>0$ satisfying
\begin{eqnarray}\label{p20120}
(1+\varepsilon)\big(a+(2+A_1)\varepsilon\big)<M,
\end{eqnarray}
there exists a constant $K=K(V,G,\gamma,\varepsilon)>0$ such that, for all $k\geq1$ and $x\ge0$,
\begin{eqnarray}\label{p201}
\overline{V^{*k}}(x)\le KM^k\overline{G}(x).
\end{eqnarray}
\end{pron}

\begin{remark}\label{remark401}
i) Clearly, $C^*(G)\ge2M(G,\gamma)$ and $M(V,\gamma)>1$, thus $a>1$. If $C^*(G)=2M(G,\gamma)$, that is $G\in\mathcal{S}(\gamma)$, then we only require $M>M(V,\gamma)$. This particular result is attributed to Lemma 2.1 of Yu et al. \cite{YWY2010}. 
And when $G=V\in\mathcal{OS}$, that is $A_1=1$, the result is due to Lemma 6.3 (ii) of Watanabe \cite{W2008}. In the two results, the distribution $V$ is supported on $(-\infty,\infty)$.

ii) Clearly, in Theorem \ref{thm1}, if $\sum_{k=1}^\infty p_kM^{k-1}<\infty,$
then by Proposition \ref{p2} with $G=V^{*\tau}$, condition (\ref{thm103}) is satisfied. For example, we can take $p_k=pq^{k},k\ge0$, where $p,q>0$ and $p+q=1$, if $q$ is small enough such that
$$qa=q\big(M(V,\gamma)+A_1(C^*(V^{*\tau})-2M(V^{*\tau},\gamma))\big)<1,$$
then we can choose a M that is small enough to make (\ref{p201}) satisfied. Thus, (\ref{thm103}) holds.

iii) From (\ref{p20110}), 
we know that $M(V,\gamma)<M-b$. Thus for all integer $n\ge1$ and any constant $K>0$,
\begin{eqnarray}\label{p20130}
M^n(V,\gamma)<M^n\big((M-b)/M\big)^n\le KM^n\big((M-b)/(KM)\big).
\end{eqnarray}
\end{remark}
\proof  Clearly, (\ref{p201}) holds for $k=1$ and all $x\ge0$. Further, we assume that (\ref{p201}) holds for $k=n$ and all $x\ge0$. For the above mentioned $\varepsilon>0$ and any $h\in{\cal{H}}(G,\gamma)$, by $G\in\mathcal{L}(\gamma)\cap\mathcal{OS}$,
there is a constant $x_0>0$ such that, for all $x\ge x_0$,
\begin{eqnarray}\label{p5}
\int_{0}^{h(x)}\overline{G}(x-y)V^{*n}(dy)\le (1+\varepsilon)M^n(V,\gamma)\overline{G}(x)\ \text{uniformly for
all}\ n\ge1,
\end{eqnarray}
\begin{eqnarray}\label{p6}
\int_{h(x)}^{x-h(x)}\overline{G}(x-y)G(dy)\leq (1+\varepsilon)\big(C^*(G)-2M(G,\gamma)+\varepsilon\big)\overline{G}(x)
\end{eqnarray}
and
\begin{eqnarray}\label{p7}
\overline{V}\big(h(x)\big)\overline{G}\big(x-h(x)\big)<\varepsilon\overline G(x).
\end{eqnarray}

For the $\varepsilon>0$, we take
$$K\ge\max\{A_1(M-b)/(M\varepsilon),1/\overline{G}(x_0)\}.$$
Further, by (\ref{p20130}), we have
\begin{eqnarray}\label{p20140}
A_1M^n(V,\gamma)\le \varepsilon KM^n\ \ \text{for all}\ n\ge1.
\end{eqnarray}

Now, we prove that (\ref{p201}) holds for $k=n+1$ and all $x\ge0$.

For all $x\ge x_0$, using integration by parts and inductive hypothesis, by (\ref{p20110}), 
(\ref{p5})-(\ref{p20140}), we have
\begin{eqnarray*}
&&\overline{V^{*(n+1)}}(x)=\int_{0}^{h(x)}\overline{V}(x-y)V^{*n}(dy)
+\int_{0}^{x-h(x)}\overline{V^{*n}}(x-y)V(dy)
+\overline{V^{*n}}\big(h(x)\big)\overline{V}\big(x-h(x)\big)\nonumber\\
&\leq&A_1\int_{0}^{h(x)}\overline{G}(x-y)V^{*n}(dy)+A_n\Big(\int_{0}^{x-h(x)}\overline{G}(x-y)V(dy)
+\overline{G}\big(h(x)\big)\overline{V}\big(x-h(x)\big)\Big)\nonumber\\
&=&A_1\int_{0}^{h(x)}\overline{G}(x-y)V^{*n}(dy)+A_n\Big(\int_{0}^{h(x)}\overline{G}(x-y)V(dy)\nonumber\\
&&\ \ \ \ \ \ \ \ \  \ \ \ \ \ \ \ \ \ \ \ +\int_{h(x)}^{x-h(x)}\overline{V}(x-y)G(dy)+\overline{V}\big(h(x)\big)\overline{G}\big(x-h(x)\big)\Big)\nonumber\\
&\leq&A_1(1+\varepsilon)M^n(V,\gamma)\overline{G}(x)
+A_n\Big((1+\varepsilon)M(V,\gamma)\overline{G}(x)+A_1\int_{h(x)}^{x-h(x)}\overline{G}(x-y)G(dy)
+\varepsilon\overline G(x)\Big)\nonumber\\
&\leq&(1+\varepsilon)\overline{G}(x)KM^n\Big(M(V,\gamma)
+A_1\big(C^*(G)-2M(G,\gamma)\big)+(2+A_1)\varepsilon\Big).
\end{eqnarray*}
Further, by (\ref{p20120}), we know that
\begin{eqnarray*}
&&A_{n+1}\le KM^{n}(1+\varepsilon)(a+2\varepsilon)\le KM^{n+1}.
\end{eqnarray*}

And for all $0\le x\le x_0$, by $M>1$ and $K\ge1/\overline{G}(x_0)$, we have
\begin{eqnarray*}
&&A_{n+1}\le 1/\overline{G}(x_0)\le KM^{n+1}.
\end{eqnarray*}

Combined with the above two inequalities, we immediately know that inequality (\ref{p201}) holds for $k=n+1$ and all $x\ge0$. \hfill$\Box$\\

\vspace{0.2cm}

\end{document}